\documentclass{amsart}

\usepackage{microtype}
\usepackage{amssymb}
\usepackage{mathtools}
\usepackage{tikz-cd}
\usepackage{mathbbol} 

\usepackage[
	backend=biber,
	style=alphabetic,
	backref=true,
	url=false,
	doi=false,
	isbn=false,
	eprint=false]{biblatex}

\renewbibmacro{in:}{}  
\AtEveryBibitem{\clearfield{pages}}  

\DeclareFieldFormat{title}{\myhref{\mkbibemph{#1}}}
\DeclareFieldFormat
[article,inbook,incollection,inproceedings,patent,thesis,unpublished]
{title}{\myhref{\mkbibquote{#1\isdot}}}

\newcommand{\doiorurl}{%
	\iffieldundef{url}
	{\iffieldundef{eprint}
		{}
		{http://arxiv.org/abs/\strfield{eprint}}}
	{\strfield{url}}%
}

\newcommand{\myhref}[1]{%
	\ifboolexpr{%
		test {\ifhyperref}
		and
		not test {\iftoggle{bbx:eprint}}
		and
		not test {\iftoggle{bbx:url}}
	}
	{\href{\doiorurl}{#1}}
	{#1}%
}

\usepackage[
	bookmarks=true,
	linktocpage=true,
	bookmarksnumbered=true,
	breaklinks=true,
	pdfstartview=FitH,
	hyperfigures=false,
	plainpages=false,
	naturalnames=true,
	colorlinks=true,
	pagebackref=false,
	pdfpagelabels]{hyperref}

\hypersetup{
	colorlinks,
	citecolor=blue,
	filecolor=blue,
	linkcolor=blue,
	urlcolor=blue
}

\usepackage[capitalize, noabbrev]{cleveref}
\crefname{subsection}{\S\!}{subsections}

\calclayout

\makeatletter
\@namedef{subjclassname@2020}{%
	\textup{2020} Mathematics Subject Classification}
\makeatother

\theoremstyle{definition}

\newtheorem*{definition*}{Definition}

\newtheorem*{remark*}{Remark}

\newtheorem*{example*}{Example}

\newtheorem*{convention*}{Convention}

\newtheorem*{notation*}{Notation}

\newtheorem*{question*}{Question}

\DeclareMathOperator{\bd}{\partial}
\DeclareMathOperator{\sign}{sign}
\newcommand{\ot}{\otimes}

\newcommand{\Z}{\mathbb{Z}}
\newcommand{\Q}{\mathbb{Q}}
\newcommand{\R}{\mathbb{R}}
\renewcommand{\k}{\Bbbk}
\newcommand{\Sym}{\mathbb{S}}
\newcommand{\Cyc}{\mathbb{C}}
\newcommand{\Ftwo}{{\mathbb{F}_2}}
\newcommand{\Fp}{{\mathbb{F}_p}}

\newcommand{\gcube}{\mathbb{I}}

\newcommand{\simplex}{\triangle}
\newcommand{\sSet}{\mathsf{sSet}}
\newcommand{\cube}{\square}

\DeclareMathOperator{\cobar}{\mathbf{\Omega}}
\DeclareMathOperator{\proj}{\pi}
\DeclareMathOperator{\incl}{\iota}

\DeclareMathOperator{\chains}{N}

\DeclareMathOperator{\gchains}{C}

\DeclarePairedDelimiter\bars{\lvert}{\rvert}

\newcommand{\id}{\mathsf{id}}
\renewcommand{\th}{\mathrm{th}}
\newcommand{\op}{\mathrm{op}}
\DeclareMathOperator*{\colim}{colim}

\DeclareMathOperator{\Med}{\mathcal{M}}
\newcommand{\Hom}{\mathrm{Hom}}
\newcommand{\End}{\mathrm{End}}

\newcommand{\xra}[1]{\xrightarrow{#1}}

\newcommand{\bC}{\mathbb{C}}

\newcommand{\bQ}{\mathbb{Q}}

\newcommand{\cL}{\mathcal{L}}

\newcommand{\cW}{\mathcal{W}}

\newcommand{\rH}{\mathrm{H}}

\newcommand{\rP}{\mathrm{P}}

\newcommand{\pdfEinfty}{\texorpdfstring{E_\infty}{E-infty}}
\usepackage{csquotes}

\DeclareMathOperator{\UMed}{\mathrm{U}(\mathcal{M})}
\DeclareMathOperator{\copr}{\Delta}
\DeclareMathOperator{\aug}{\epsilon}
\DeclareMathOperator{\pr}{\ast}
\DeclareMathOperator{\ad}{ad}
\DeclareMathOperator{\barconst}{\mathbf{B}}
\newcommand{\desus}[1]{s^{\mbox{\scriptsize -#1}}}
\DeclareMathSymbol{\antishrik}{\mathord}{operators}{"3C}

\newcommand{\cdga}{\mathsf{cdga}}
\newcommand{\cdgl}{\mathsf{cdgl}}
\newcommand{\com}{\mathcal{C}om}
\newcommand{\lie}{\mathcal{L}ie}
\newcommand{\Der}{\mathrm{Der}}
\title[Effective algebro-homotopical constructions]
{The diagonal of cellular spaces and effective algebro-homotopical constructions}

\author{Anibal~M.~Medina-Mardones}
\address{A.M-M., Max Planck Institute for Mathematics \and University of Notre Dame}
\email{\href{mailto:ammedmar@mpim-bonn.mpg.de}{ammedmar@mpim-bonn.mpg.de}}

\dedicatory{In honor of Dennis Sullivan's $80^\th$ birthday}
\date{\today}
\subjclass[2020]{55U15, 18N40, 18M70, 18M85, 55S30, 55S05}
\keywords{Cellular chains, diagonal map, effective constructions, $C_\infty$-coalgebras, $E_\infty$-coalgebras}

\begin{document}
	\begin{abstract}

In this survey article we discuss certain homotopy coherent enhancements of the coalgebra structure on cellular chains defined by an approximation to the diagonal.
Over the rational numbers, $C_\infty$-coalgebra structures control the $\bQ$-complete homotopy theory of spaces, and over the integers, $E_\infty$-coalgebras provide an appropriate setting to model the full homotopy category.
Effective constructions of these structures, the focus of this work, carry geometric and combinatorial information which has found applications in various fields including deformation theory, higher category theory, and condensed matter physics.
	\end{abstract}
	\maketitle

\section{Introduction} \label{s:introduction}

There is a tense trade-off in algebraic topology having roots reaching back to the beginning of its modern form.
This tension can be illustrated with the concept of cohomology.
The first approaches, dating back to Poincar\'e, are based on the subdivision of a space into simple contractible pieces.
These elementary shapes are made to generate a free graded module and their spatial relations define the differential used to compute it.
This definition makes certain geometric properties of cohomology, for example excision, fairly clear.
Yet, it is not easy to show that a continuous map of spaces induces a map between their associated cohomologies.
The functoriality just alluded to is trivial when defining cohomology in terms of homotopy classes of maps to Eilenberg--MacLane spaces, but the passage to the homotopy category erases geometric and combinatorial information and the resulting definition is not well suited for concretely presented spaces.

Cohomology as a graded group is a fairly computable invariant but it has noticeable limitations, for example, $\bC \rP^2$ and $S^2 \vee S^4$ are not distinguished by it.
This invariant can be refined to a graded ring by endowing cohomology with the cup product.
In the spectral context, this finer structure is defined through the wedge product of Eilenberg--MacLane spaces, whereas in the cellular setting it is obtained from a choice of
cellular approximation to the diagonal.
Such cellular map induces a chain map
\begin{equation} \label{e:chain diagonal}
\Delta \colon \gchains(X) \to \gchains(X)^{\otimes 2}
\end{equation}
making the cellular chains of $X$ into a (differential graded) coalgebra.
The fact that the cup product on cohomology, induced by the linear dual of \eqref{e:chain diagonal}, is associative and (graded) commutative hints at the presence of additional structure extending the coalgebra structure on $\gchains(X)$.

In this survey article we will present, from a viewpoint that emphasizes their constructive nature, $C_\infty$ and $E_\infty$ extensions of $\Delta$ over the rationals and integers respectively.
The resulting algebraic structures control much of the homotopy theory of spaces.
For example, over the rationals, the quasi-isomorphism type of a $C_\infty$-coalgebra extension of the symmetrization of $\Delta$ determines the $\Q$-completion of $X$ under certain assumptions \cite{quillen1969rational, buijs2020liemodels}.
Similarly, over the algebraic closure of $\Fp$, the quasi-isomorphism type of the dual of an $E_\infty$-coalgebra extension of $\Delta$ determines the $p$-completion of $X$ provided certain finiteness assumptions are met \cite{mandell2001padic}.
Both of these statements follow from the stronger fact that the relevant model categories are Quillen equivalent under certain assumptions.

\subsection*{Rational coefficients}

In \cref{s:rationally} we will study extensions of chain approximations to the diagonal with rational coefficients.
Over this field, a chain approximation to the diagonal can be symmetrized, giving rise to a cocommutative coalgebra.
This coalgebra cannot be made simultaneously coassociative, but this relation can be imposed in a derive sense through a family of coherent chain homotopies -- which also respect certain symmetry constrains -- and give rise to a so called $C_\infty$-coalgebra structure.
One can think of $C_\infty$-coalgebras in terms of the somewhat more familiar notion of $A_\infty$-coalgebra where cocommutativity is satisfied strictly.
As a manifestation of Koszul duality, a $C_\infty$-coalgebra structure on cellular chains is equivalent to a differential on the completion of the free graded Lie algebra generated by the cells shifted downwards in degree by one.
This relates $C_\infty$-coalgebras to deformation theory, but we do not explore this deep connection here.
For cell complexes whose closed cells have the $\Q$-homology of a point, Dennis provided in \cite{sullivan2007appendix} a local inductive construction defining a $C_\infty$-coalgebra structure on their cellular chains.
We reprint a challenge he posted regarding the resulting structure.
\begin{displaycquote}[p.2]{lawrence2014interval}
	\textsc{Problem}. Study this free differential Lie algebra attached to a cell complex, determine its topological and geometric meaning as an intrinsic object.
	Give closed form formulae for the differential and for the induced maps associated to subdivisions.
\end{displaycquote}
As proven by Quillen, the quasi-isomorphism type of this $C_\infty$-coalgebra is a complete invariant of the rational homotopy type of simply-connected spaces.
For the $C_\infty$-coalgebra structure on the interval, Dennis and Ruth Lawrence addressed the challenge reprinted above introducing a formula for it which can be interpreted in terms of parallel transport of flat connections \cite{lawrence2014interval}, and for which the subdivision map is described by the Baker--Campbell--Hausdorff formula.

To generalize Quillen's equivalence of homotopy categories to one between (not necessarily 1-connected) simplicial sets, Buijs, F{\'e}lix, Murillo, and Tanr{\'e} \cite{buijs2020liemodels} extended to $n$-simplices the Lawrence--Sullivan structure building, constructively for $n \in \{2,3\}$ and inductively otherwise, $C_\infty$-coalgebra structures on their chains.
Their construction agrees after linear dualization with the one obtained by Cheng and Getzler in \cite{getzler2008transfering}, where they showed that the Kontsevich--Soibelman sum-over-trees formula defining the transfer of $A_\infty$-algebras through a chain contraction induces a transfer of $C_\infty$-algebras.
This allowed them to construct a $C_\infty$-algebra structure on simplicial cochains by transferring Dennis' polynomial differential forms through Dupont's contraction.
The resulting description is given in terms of rooted trees.

$C_\infty$-coalgebras are controlled by the operad $\com_\infty$ which is the Koszul resolution of the operad $\com$, i.e., the cobar construction applied to the $\lie$ cooperad, the Koszul dual cooperad of $\com$.
Another interesting resolution of $\com$ is constructed concatenating the bar and cobar constructions.
This resolution method is an algebraic version of the $W$-construction of Boardman--Vogt.
As Dennis and Scott Wilson considered, the resulting operad can be described using rooted trees with vertices colored black or white.
In \cite{vallette2020higherlietheory}, Daniel Robert-Nicoud and Bruno Vallette studied coalgebras over this resolution which they termed $CC_\infty$-coalgebras.
They constructed on the chain of standard simplices natural $CC_\infty$-coalgebra structures and described them explicitly using bicolored trees.

Despite some progress -- \cite{lawrence2019triangle, lawrence2021cells, buijs2019triangle}, \cite[\S6.5]{buijs2020liemodels} -- the ``closed form formulae'' part of the problem quoted before remains open.
One possible avenue to generalize to cubical chains the formula defining the Lawrence--Sullivan $C_\infty$-coalgebra on $\gchains(\gcube)$, is to define the tensor product of $C_\infty$-coalgebras and then extend it monoidally to all cubes via the isomorphism $\chains(\cube^n) \cong \gchains(\gcube)^{\ot n}$.
The monoidal structure on the category of $A_\infty$-coalgebras is defined through a chain approximation to the diagonal of the Stasheff polytopes compatible with the operad structure.
Unfortunately, the resulting $A_\infty$-coalgebra on $\gchains(\gcube)^{\ot 2}$ is not $C_\infty$.
This could be corrected through an algebraic symmetrization of the associahedral diagonal, but we do not pursue this here.

\subsection*{Integral coefficients}

In \cref{s:integrally} we will study extensions of chain approximations to the diagonal with integral coefficients.
In contrast to the situation over $\Q$, chain approximations to the diagonal over these coefficients cannot be taken to be symmetric with respect to transposition of tensor factors.
The resulting coalgebras can be made cocommutative and coassociative only up to coherent homotopies, that is to say, provided with the structure of a coalgebra over an $E_\infty$-operad.
The study of $E_\infty$-structures has a long history, where (co)homology operations \cite{steenrod1962cohomology, may1970general}, the recognition of infinite loop spaces \cite{boardman1973homotopy, may1972geometry}, and the complete algebraic representation of the $p$-adic homotopy category \cite{mandell2001padic} are key milestones.

Steenrod was the first to introduce homotopy coherent corrections to the broken symmetry of a chain approximation to the diagonal \cite{steenrod1947products}.
He did so on simplicial chains in the form of explicit formulae defining his cup-$i$ coproducts, with \mbox{cup-$0$} agreeing with the Alexander--Whitney chain approximation to the diagonal.
These coproducts are used to define Steenrod's mod $2$ cohomology operations and to effectively compute them in specific examples.

Extending the cup-$i$ coproducts of Steenrod, explicitly defined $E_\infty$-coalgebra structure on simplicial chains were introduced by McClure--Smith \cite{mcclure2003multivariable} and Berger--Fresse \cite{berger2004combinatorial}.
It turns out that this structure can be describe solely in terms of the Alexander--Whitney diagonal, the augmentation map and a chain version of the join of simplices \cite{medina2020prop1}.
This point of view can be abstracted using the language of props, which allows its application to other contexts, for example those defined by cubical chains \cite{medina2021cubical} and the Adams' cobar construction \cite{medina2021cobar}.
We will review the resulting model of the $E_\infty$-operad, its action on simplicial and cubical chains, and explicit generalizations of the cup-$i$ coproducts to higher arities effectively constructing Steenrod operations at all primes \cite{medina2021may_st}.

We devote the final subsection to one use of cochain level structures in the classification of symmetry protected topological phases of matter.


\subsection*{Acknowledgment}

We would like to thank Greg Brumfiel, Ruth Lawrence, John Morgan, Aniceto Murillo, Dennis Sullivan and Bruno Vallette for helpful exchanges related to their work.
We are also grateful to the referees for their careful reading and many keen suggestions.

\section{\texorpdfstring{$C_\infty$}{C-infty}-coalgebras} \label{s:rationally}

Over the rationals, the problem of extending a chain approximation to the diagonal as a $C_\infty$-coalgebra is related to the study of Lie algebras.
In this section we recall this connection, and Dennis' construction of a $C_\infty$-coalgebra structure on the cellular chains of certain CW complexes.
We also discuss the resulting structure on the cellular chains of the interval, which is presented as a formula in the work of Dennis and Ruth Lawrence.
We discuss Quillen's functor from simplicial sets to complete dg Lie algebras, as extended by Buijs, F{\'e}lix, Murillo, and Tanr{\'e} through a cosimplicial $C_\infty$-coalgebra, and the problem of making this construction into explicit formulae extending the Sullivan--Lawrence interval.

\subsection{Quillen construction}

As a motivating example illustrating the connection between cocommutative and coassociative coalgebras and dg Lie algebras, let us recall the so-called \textit{Quillen construction}.
Consider one such coalgebra $C$, and form the free graded Lie algebra $L$ generated by the desuspension of $C$ regarded as a graded vector space.
Denote $\desus{1}C$ in $L$ by $L_1$, and by $L_2$ the linear span of brackets of elements in $\desus{1}C$.
The boundary map and coproduct respectively induce maps
\[
L_1 \xra{l_1} L_1,
\qquad
L_1 \xra{l_2} L_2,
\]
of degree $-1$, and their relations ensure that $l_1 + l_2$ squares to $0$.
More explicitly,
\[
l_1(\desus{1}c) = -\desus{1} \bd c,
\qquad
l_2(\desus{1}c) = \frac{1}{2} \sum_i (-1)^{\bars{a_i}} \big[ \desus{1}a_i, \desus{1}b_i \big],
\]
where $\Delta(c) = \sum_i a_i \otimes b_i$.

The extension of $l_1 + l_2$ as a derivation of the Lie bracket makes $L$ into a free dg Lie algebra naturally associated to $C$.

\subsection{\texorpdfstring{$C_\infty$}{C-infinity}-coalgebras} \label{ss:c-infty definition}

The previous construction motivates the definition of $C_\infty$-coalgebras.
Before providing it, let us recall the notion of a \textit{complete chain complex} $(C, F)$, which is a filtered chain complex $C$
\[
C = F_0 C \supseteq F_1 C \supseteq \cdots
\]
such that
\[
C = \lim_{k\to \infty} C / F_k C.
\]
As expected, the \textit{completion} of a filtered chain complex $(C, F)$ is defined as $\lim_{k\to \infty} C / F_k C$.

A $C_\infty$-\textit{coalgebra structure} on a graded vector space $C$ is the data of a differential on the completion, with respect to the filtration by number of brackets, of the free graded Lie algebra generated by $\desus{1}C$.

\subsection{$C_\infty$-coalgebras as commutative $A_\infty$-algebras} \label{ss:a-infty coalgebras}

We can interpret a $C_\infty$-coalgebra structure on $C$ in terms of the somewhat more familiar notion of $A_\infty$-coalgebra.

An $A_\infty$-\textit{coalgebra structure} on a graded vector space $C$ is a family of degree $k-2$ linear maps $\Delta_k \colon C \to C^{\ot k}$ satisfying for every $i \geq 1$ the following identity:
\begin{equation} \label{e:a-infty relations}
\sum_{k=1}^{i} \sum_{n=0}^{i-k} (-1)^{k+n+kn} \big( \id^{\ot i-k-n} \ot \Delta_k \ot \id^{\ot n} \big) \circ \Delta_{i-k+1} = 0.
\end{equation}
This is equivalent to the data of a differential on $\prod_{n \geq 1} (\desus{1} C)^{\otimes n}$, the augmentation kernel of the complete tensor algebra on the desuspension of $C$.
Indeed, such differential $d = \sum_{k \geq 1} d_k$ is determined by its restriction to $\desus{1} C$ with $d_k(\desus{1} C) \subset T^k(\desus{1} C)$, and the correspondence is explicitly given by
\[
\Delta_k = - s^{\ot k} \circ d_k \circ \desus{1},
\qquad
d_k = -(-1)^{\frac{k(k+1)}{2}} (\desus{1})^{\ot k} \circ \Delta_k \circ s.
\]
Notice that \eqref{e:a-infty relations} implies for any $A_\infty$-coalgebra that $\Delta_1$ squares to $0$, that $\Delta_2$ is a chain map with respect to $\Delta_1$, and that $\Delta_3$ is a chain homotopy between $(\Delta_2 \otimes \id) \circ \Delta_2$ and $(\id \ot \Delta_2) \circ \Delta_2$.

A $C_\infty$-coalgebra structure on a graded vector space $C$ is equivalent to an $A_\infty$-coalgebra structure on $C$ such that the image $d_k$ lies in the invariants of $(\desus{1} C)^{\ot k}$ under the action of $\Sym_k$, or, expressed in terms of the coproducts $\Delta_k$, one such that $\tau \circ \Delta_k = 0$, where
\[
\tau(c_1 \ot \dotsb \ot c_k) =
\sum_{i=1}^{k} \sum_{\sigma \in \Sym(i, k-i)}
\!\!\! \sign(\sigma) \,
(c_{\sigma(1)} \ot \dotsb \ot c_{\sigma(i)}) \ot
(c_{\sigma(i+1)} \ot \dotsb \ot c_{\sigma(k)})
\]
and $\Sym(i, k-i)$ denotes the set of $(i, k-i)$-shuffles.

\subsection{Sullivan's cellular $C_\infty$-coalgebra construction} \label{ss:dennis construction}

We now present Dennis' inductive construction of a local $C_\infty$-coalgebra structure on the chains of cell complexes whose closed cells have the $\Q$-homology of a point \cite{sullivan2007appendix}.

Let $X$ be one such cellular complex and $L(X) = L$ be the free Lie algebra generated by the desuspension of its rational cellular chains $\desus{1}C$.

Let us start by choosing a chain approximation $\copr \colon C \to C \ot C$ to the diagonal, which we assume equivariant -- since we are working with rational coefficient -- and local, in the sense that $\copr(e_\alpha)$ is contained in the subcomplex generated by the tensor product of cells in the closure of $e_\alpha$.
We remark that $(C, \delta, \copr)$ is a cocommutative coalgebra which is in general not coassociative.
Let $\delta_1$ and $\delta_2$ be the respective maps from $L_1$ to $L_1$ and $L_2$ induced from $\bd$ and $\copr$.
We denote by the same symbols their extensions to $L$ as derivations.
We now quote Dennis' construction:

\begin{displaycquote}[p.251]{sullivan2007appendix}
	Interpreting the equation $\delta \circ \delta = 0$ as $[\delta, \delta] = 0$ where $[\cdot, \cdot]$ is the graded commutator.
	For any $\delta$ the Jacobi identity is $[\delta, [\delta, \delta]]$.
	Suppose $\delta^k = \delta_1 + \dots + \delta_k$ has been defined so that $[\delta^k, \delta^k]$ has the first nonzero term in monomial degree $k + 1$.
	Jacobi implies this error commutes with $\delta_1$; that is, it is a closed element in the complex $\Der(L)$ of derivations of $L$.
	If we work in the closure of a cell, the homology hypothesis implies that $\Der(L)$ has homology only in degrees $0$ and $1$.
	Therefore, the error, which lives in degree $2$, can be written as a commutator with $\delta_1$.
	Using the cells to generate a linear basis of each $L_k$ by bracketing, we choose this solution to lie in the image of the adjoint of $\delta_1$ to make it canonical.
	This canonical solution is $\delta_{k+1}$ and this completes the induction, since one knows at the beginning $\delta_1 \circ \delta_1$
	and $\delta_2$ is chain mapping; that is, $[\delta_2, \delta_1] = 0$.
\end{displaycquote}

Dennis' construction is such that $\delta e_\alpha$ is in the sub Lie algebra generated by the closure of the cells in $e_\alpha$, or, expressed in dual terms, the maps $\copr_r \colon C \to C^{\ot r}$ corresponding to the $\delta_r$ maps are local.

Recall Dennis' problem, quoted in the introduction, of determining the topological and geometric meaning of this $C_\infty$-coalgebra as an intrinsic object, and give closed form formulae for it and the induced maps associated to subdivisions.
We will next present the solution Dennis and Ruth Lawrence gave to this problem in the case of the interval.

\subsection{Lawrence--Sullivan interval} \label{ss:LS interval}

Let $L$ be a completed free graded Lie algebra with filtration given by number of brackets, and let $U(L)$ be the complete graded vector space of series on one indeterminate with values on $L$ whose filtration is induced from that of $L$, i.e., the $N^\th$-part of the filtration $U(L)$ contains series of the form
\[
\sum_{n=1}^\infty x_n t^n
\]
where $x_n$ is in $F_N L$ for every $n$.
Consider the linear operator given by
\[
\frac{d}{dt} \left(\sum x_n t^n\right) = \sum n \, x_n t^{n-1}
\]
and the formal differential equation
\[
\frac{du}{dt} = \bd v - \ad_v u
\]
where $\ad_v u = [v, u]$.
By formally solving this equation one defines the \textit{flow generated} by $v$ for any rational time $t_0$.

An element $u \in L$ is said to be \textit{flat} if it is in degree $-1$ and satisfies $\bd u = \frac{1}{2} [u,u]$.
It is common to refer to these as \textit{Maurer--Cartan} elements, but do not use this terminology.
We now quote Dennis and Ruth Lawrence theorem.
\begin{displaycquote}[Theorem 1]{lawrence2014interval}
	There is a unique completed free differential graded Lie algebra, $A$, with generating elements $a$, $b$ and $e$, in degrees -1, -1 and 0 respectively, for which $a$ and $b$ are flat while the flow generated by $e$ moves from $a$ to $b$ in unit time.
	The differential is specified by
	\[
	\bd e = \ad_e b + \sum_{i=0}^{\infty} \frac{B_i}{i!}(\ad_e)^i(b-a),
	\]
	where $B_i$ denotes the $i$th Bernoulli number defined as coefficients in the expansion
	\[
	\frac{x}{e^x-1} = \sum_{n=0}^{\infty} B_n \frac{x^n}{n!}.
	\]
\end{displaycquote}
We remark that Dennis conjectured an equivalence between the description above and the one obtained by applying his inductive procedure (\cref{ss:dennis construction}).
This conjecture was verified by Parent and Tanr\'{e} \cite{parent2012interval}.

The Lawrence--Sullivan dg Lie algebra is describe in terms of the associated $C_\infty$-coalgebra by
\begin{align*}
& \Delta_1(c) = y-z, \quad
\Delta_1(y) = \Delta_1(z) = 0, \\
& \Delta_2(c) = -\frac{1}{2} \Big( c \ot (y+z) + (y+z) \ot c \Big), \quad
\Delta_2(y) = -y \ot y, \quad
\Delta_2(z) = -z \ot z, \\
& \Delta_k(c) = \sum_{p+q=k-1} \frac{B_{k-1}}{p!q!} c^{\ot p} \ot (y-z) \ot c^{\ot q}, \quad
\Delta_k(y) = \Delta_k(z) = 0, \quad k \geq 3,
\end{align*}
where $c$, $y$ and $x$ are generators of degree 1, 0 and 0, respectively.

\subsection{Rational homotopy theory} \label{ss:cdgl model}

To algebraically model the rational homotopy category of spaces two models were introduced.
On one hand there is Dennis' commutative approach \cite{sullivan1977infinitesimal} based on an adjunction
\[
\begin{tikzcd}
\sSet \arrow[r, shift left=2pt, "A_{\mathrm{PL}}"] &
\cdga^\op \arrow[l, shift left=2pt, "\ \bars{\,\cdot\,}_S"].
\end{tikzcd}
\]
On the other, there is Quillen's Lie approach, introduced in \cite{quillen1969rational} and extended in \cite{buijs2013algebraicmodels, buijs2020liemodels}, which is based on an adjunction
\[
\begin{tikzcd}
\sSet \arrow[r, shift left=2pt, "\cL"] &
\cdgl \arrow[l, shift left=2pt, "\ \bars{\,\cdot\,}_Q"]
\end{tikzcd}
\]
where $\cdgl$ denotes the category of complete dg Lie algebras.
This adjunction is defined explicitly by
\[
\cL(X) = \colim_{\simplex^n \to X} \cL(\simplex^n), \qquad
\bars{L}_n = \cdgl(\cL(\simplex^n), L),
\]
where $\cL(\simplex^\bullet)$ is the cosimplicial complete dg Lie algebra defined by a natural $C_\infty$-coalgebra structures on the cellular chains of standard simplices.

Using the principles presented in the previous subsection and a careful treatment of the simplicial structure, Buijs, F{\'e}lix, Murillo, and Tanr{\'e} \cite{buijs2020liemodels} introduced a construction of $\cL(\simplex^\bullet)$, which they characterize axiomatically by requiring that the generators associated to vertices are flat, and that the linear part is induced from the boundary of chains.

We mention that this structure is isomorphic to the one obtained by dualizing the simplicial $C_\infty$-algebra defined by the Homotopy Transfer Theorem of $C_\infty$-algebras applied to Dennis' polynomial differential forms and Dupont's contraction \cite{getzler2008transfering}.

The problem of finding closed formulae for the $C_\infty$-coalgebra structure on the $n$-simplex remains open for $n > 3$.

\subsection{Operadic viewpoint}

The operad $C_\infty$ is defined as the cobar construction on the Lie cooperad, the Kozul dual cooperad of $\com$.
That is to say
\[
C_\infty = \cobar \lie^{\mathrm c}.
\]
The operad $C_\infty$ is a minimal projective resolution of $\com$.
A larger projective resolution is defined composing the bar and cobar constructions
\[
CC_\infty = \cobar \barconst \com.
\]
Since the bar and cobar constructions are defined in terms of rooted trees the $CC_\infty$ operad can be described using bicolored trees.
In \cite{vallette2020higherlietheory}, Robert-Nicoud and Vallette constructed a cosimplicial $CC_\infty$-coalgebra in terms of bicolored trees, and explored the induced adjuntion between simplicial sets and $L_\infty$-algebras.

An anecdote shared with the author by both Dennis and Bruno Vallete, is that this bicolored model brought them in contact for the first time; after a talk where Dennis used this pictorial description, Bruno, then a recent graduate, recognized it as the bar-cobar resolution of $\com$.

\section{$\pdfEinfty$-coalgebras} \label{s:integrally}

In this section we consider extensions of the Alexander--Whitney and Serre coalgebras on simplicial and cubical chains respectively.
We recall Steenrod's cup-$i$ coproducts and their induced square operations on mod 2 cohomology.
We then describe an extension of these coproducts as a full $E_\infty$-coalgebra using three linear maps, and identify generalizations to cup-$(p,i)$ coproducts inducing Steenrod operations on mod $p$ cohomology.
We close this section overviewing the use of cochain level structures in the study of symmetry protected topological phases.

\subsection{Alexander--Whitney coalgebra} \label{ss:aw diagonal}

The first chain approximation to the diagonal was given in the simplicial context by \v{C}ech and Whitney building on independent work presented during the 1935 International Congress in Moscow by Alexander and Kolmogorov.
The original references are \cite{alexander1936ring, cech1936multiplication, whitney1938products} and a historical account is presented by Whitney in \cite[p.110]{whitney1988history}.
This chain map, referred to as the \textit{Alexander--Whitney coproduct}, is defined on elements of the canonical basis by the formula
\begin{equation} \label{e:alexander-whitney coalgebra}
\copr \big( [0,\dots,n] \big) = \sum_{i=0}^n \ [0, \dots, i] \ot [i, \dots, n].
\end{equation}
Together with the \textit{augmentation map}
\begin{equation} \label{e:augmentation map}
\aug \big( [0, \dots, n] \big) =
\begin{cases}
1 & n = 0, \\
0 & n > 0,
\end{cases}
\end{equation}
the Alexander--Whitney coproduct satisfies
\begin{gather}
\label{e:coassociativity relation}
(\copr \ot \, \id) \circ \copr = (\id \ot \copr) \circ \copr, \\
\label{e:counital relation}
(\aug \ot \, \id) \circ \copr = \id = (\id \ot \aug) \circ \copr,
\end{gather}
making the (normalized) chains $\chains(X)$ of any simplicial set $X$ into a natural coassociative counital coalgebra, referred to as the \textit{Alexander--Whitney coalgebra} of $X$.

We will use the following recursively defined notation for general coalgebras:
\begin{align*}
\copr^1 &= \copr, \\
\copr^k &= (\copr \ot \, \id) \circ \copr^{k-1}.
\end{align*}

\subsection{The join product}

\begin{figure}
	\begin{tikzpicture}[scale=.6]
\coordinate (A) at (210:2);
\coordinate (B) at (-30:2);
\coordinate (C) at (90:2);

\draw[draw=black] (A) -- (B) -- (C) -- (A);

\node[circle,fill=blue, opacity=.9, inner sep=0pt,minimum size=5pt, label=left:{0}] (a) at (A) {};
\node[circle,fill=black,inner sep=0pt,minimum size=3pt, label=right:{$1$}] (a) at (B) {};
\node[circle,fill=black,inner sep=0pt,minimum size=3pt, label=right:{$2$}] (a) at (C) {};

\node[scale=1.5] at (3.5,0.5) {$\ast$};
\end{tikzpicture}
\begin{tikzpicture}[scale=.6]
\coordinate (A) at (210:2);
\coordinate (B) at (-30:2);
\coordinate (C) at (90:2);

\draw[draw=blue,  ultra thick, draw opacity=.7] (B) -- (C);
\draw[draw=black] (C) -- (A);
\draw[draw=black] (A) -- (B);

\node[circle,fill=black,inner sep=0pt,minimum size=3pt, label=left:{$0$}] (a) at (A) {};
\node[circle,fill=black,inner sep=0pt,minimum size=3pt, label=right:{$1$}] (a) at (B) {};
\node[circle,fill=black,inner sep=0pt,minimum size=3pt, label=right:{$2$}] (a) at (C) {};

\node[scale=1.5] at (3.5,.5) {=};
\end{tikzpicture}
\begin{tikzpicture}[scale=.6]
\coordinate (A) at (210:2);
\coordinate (B) at (-30:2);
\coordinate (C) at (90:2);

\draw[draw=black] (A) -- (B) -- (C) -- (A);

\node[circle,fill=black,inner sep=0pt,minimum size=3pt, label=left:{$0$}] (a) at (A) {};
\node[circle,fill=black,inner sep=0pt,minimum size=3pt, label=right:{$1$}] (a) at (B) {};
\node[circle,fill=black,inner sep=0pt,minimum size=3pt, label=right:{$2$}] (a) at (C) {};

\draw[draw, fill=blue, opacity=.7] (A) -- (B) -- (C) -- (A);
\end{tikzpicture}
	\caption{Geometric representation of the join product of two basis elements. It depicts the identity $\pr \big( [0] \otimes [1,2] \big) = [0,1,2]$.}
	\label{f:join of faces}
\end{figure}
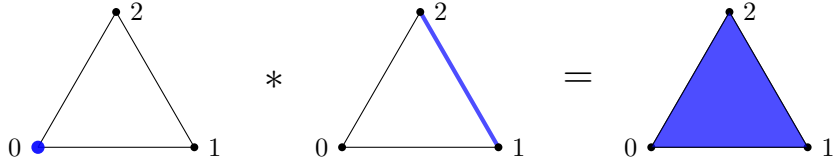

The \textit{join product} $\ast \colon \chains(\simplex^n)^{\ot 2} \to \chains(\simplex^n)$ is the natural degree~$1$ linear map defined by
\begin{multline}
\ast \big(\left[v_0, \dots, v_p \right] \ot \left[v_{p+1}, \dots, v_q\right]\big) = \\
\begin{cases} (-1)^{p} \sign(\pi) \left[v_{\pi(0)}, \dots, v_{\pi(q)}\right] & \text{ if } \forall\, i \neq j, \ v_i \neq v_j, \\
\hfil 0 & \text{ if not}, \end{cases}
\end{multline}
where $\pi$ is the permutation that orders the vertices.
It is an algebraic version of the usual join of faces in a simplex, please consult \cref{f:join of faces} for an example.

The join product can be used in conjunction with the Alexander--Whitney coproduct to canonically construct boundaries in the chain complexes
\[
\Hom \big( \chains(\simplex^n)^{\ot s}, \chains(\simplex^n)^{\ot r} \big).
\]
For example,
\[
H = (f \ast g) \circ \copr
\]
is a chain homotopy between any two quasi-isomorphisms $g, f \colon \chains(\simplex^n) \to \chains(\simplex^n)$.
To see this, recall the augmentation map $\aug \colon \chains(\simplex^n) \to \k$ defined in \eqref{e:augmentation map} which is the counit of $\copr$, and notice that the join is a chain homotopy between $\aug \ot \, \id$ and $\id \ot \aug$, that is to say
\begin{equation}
\partial \pr = \aug \ot \, \id - \id \ot \aug.
\end{equation}
Since $f$ and $g$ are quasi-isomorphisms we have $\aug \circ f = \aug \circ g = \aug$, so
\begin{align*}
\bd H &=
\big( \aug \ot \, \id - \id \ot \aug \big) \circ (f \ot g) \circ \copr \\ &=
\big(\aug \ot \, g - f \ot \aug \big) \circ \copr \\ &= g - f.
\end{align*}

\subsection{Steenrod cup-$i$ coproduct structure} \label{ss:cup-i}

As it can be seen directly from \eqref{e:alexander-whitney coalgebra}, the Alexander--Whitney coproduct is not cocommutative.
In \cite{steenrod1947products}, Steenrod introduced coherent higher diagonals correcting homologically this lack of cocommutativity.
He used them to define the celebrated square operations, finer invariants on the mod~2 cohomology of spaces (\cref{ss:steenrod squares}).
In this subsection we present an explicit recursive definition of Steenrod's higher diagonals.

Let $C$ be a chain complex of $\Z$-modules and regard $\Hom(C, C \ot C)$ as a chain complex of $\Z[\Sym_2]$-modules where $\Sym_2$ acts by permuting the factors in the target.
Denote the elements $1 + (12)$ and $(12) - 1$ in $\Z[\Sym_2]$ by $N$ and $T$ respectively.
A \textit{cup-$i$ coproduct structure} on $C$ is an equivariant chain map
$\cW(2) \to \Hom(C, C \ot C)$ where
\begin{equation} \label{e:minimal resolution r=2}
\begin{tikzcd} [column sep = .5cm]
\mathcal W(2) = \Z[\Sym_2]\{e_0\} & \arrow[l, "\,T"'] \Z[\Sym_2]\{e_1\} & \arrow[l, "\,N"'] \Z[\Sym_2]\{e_2\} & \arrow[l, "\,T"'] \cdots
\end{tikzcd}
\end{equation}
is the minimal free resolution of $\Z$ as a $\Z[\Sym_2]$-module.
The image of $e_i$ is denoted by $\copr_i \colon C \to C \otimes C$ and is referred to as the \textit{cup-$i$ coproduct} of $C$ (with respect to the given cup-$i$ coproduct structure).

We can use the Alexander--Whitney coproduct and the join product to give a recursive description of the natural cup-$i$ coproduct structure on simplicial chains introduced in \cite[p.293]{steenrod1947products}:
\begin{equation} \label{e:cup-i coproducts}
\begin{split}
& \copr_0 = \copr, \\
& \copr_i =
(\ast \ot \id) \circ (\id \ot (12)\copr_{i-1}) \circ \copr.
\end{split}
\end{equation}

We refer to \cite{mcclure2003multivariable, gonzalez-diaz1999steenrod, medina2021fast_sq} for alternative descriptions of isomorphic \mbox{cup-$i$} constructions, where we say that two \mbox{cup-$i$} constructions on $C$, say $\psi$ and $\psi^\prime$, are \textit{isomorphic} if there is an automorphism of $\cW(2)$ making the following diagram commute:
\[
\begin{tikzcd} [column sep = -25, row sep=normal]
\cW(2) \arrow[rr, "\cong"] \arrow[rd, in=150, out=-90, "\psi^{\phantom{\prime}}"', near start] & & \cW(2) \arrow[ld, in=30, out=-90, "\psi^\prime", near start] \\
& \Hom(C, C^{\ot 2}). &
\end{tikzcd}
\]
The cup-$i$ products of Steenrod seem to be combinatorially fundamental, for example, they induce the nerve of strict infinity categories \cite{street1987orientals, medina2020globular} and can be axiomatically characterized \cite{medina2022axiomatic}.

\subsection{Steenrod square operations} \label{ss:steenrod squares}

Let $C$ be equipped with a cup-$i$ coproduct structure.
The \textit{Steenrod square operations}
\[
Sq^k \colon \rH(C^\vee) \to \rH(C^\vee)
\]
on the homology of its dual chain complex $C^\vee = \Hom(C, \Ftwo)$ are defined for every integer $k$ by the formula
\begin{equation} \label{e:steenrod squares}
Sq^k \big( [\alpha] \big) = \big[ (\alpha \ot \alpha) \copr_{k - \bars{\alpha}}(-) \big]
\end{equation}
where brackets are used to denote represented elements in $\rH(C^\vee)$.

\subsection{An $\pdfEinfty$-coalgebra on simplicial chains} \label{ss:e-infty generalization}

Cup-$i$ coproducts on simplicial chains are part of an $E_\infty$-coalgebra structure.
This is a natural coalgebra structure over an operad whose arity $r$ part is a chain complex of free $\k[\Sym_r]$-module with the $\k$-homology of a point.
Similar to Dennis' construction over $\Q$ of an $C_\infty$-coalgebra structure on cellular chains (\cref{ss:dennis construction}), the existence of an $E_\infty$-coalgebra structure over any coefficient ring can be guaranteed using an acyclic carrier argument \cite{eilenberg1953acyclic}.
The goal of this subsection is to describe explicitly an $E_\infty$-coalgebra structure on simplicial integral chains generalizing the construction of cup-$i$ coproducts of Steenrod (\cref{ss:cup-i}).

The collection of all linear maps $\chains(\simplex^n) \to \chains(\simplex^n)^{\ot r}$ for any $r$ that can be expressed as an arbitrary compositions of the Alexander--Whitney coproduct, the join product, and permutations of factors defines an $E_\infty$-coalgebra structure on the chains of standard simplices.
We remark that, since we are only considering maps whose domain is $\chains(\simplex^n)$, the join is not part of this structure, although it is used in its construction.

The $E_\infty$-operad $\UMed$ defining this structure can be abstracted from this example.
Roughly speaking, $\UMed = \{\Med(1,r)\}_{r \geq 0}$ is the operad associated to the prop $\Med$ generated by symbols $\copr, \aug, \pr$ in biarities $(1,2)$, $(1,0)$, and $(2,1)$ of degree $0,0,1$ with $\bd \copr = 0$, $\bd \aug = 0$, and $\bd \pr = \aug \ot \, \id - \id \ot \aug$, modulo the relations $(\aug \ot \, \id) \circ \copr = \id = (\id \ot \aug) \circ \copr$ and $\aug \circ \, \ast = 0$.
In \cref{ss:homology of M} we review a family of explicit chain contractions that can be used to compute the homology of $\UMed$.
We use this family in \cref{ss:higher cup-i coproducts} to define cup-$(r,i)$ coproducts responsible for Steenrod operations at all primes.

Full details regarding the construction of the prop $\UMed$ can be found in \cite{medina2020prop1, medina2021prop2} together with a comparison to the surjection operad \cite{mcclure2003multivariable, berger2004combinatorial}, a construction based on an earlier generalization of Steenrod's cup-$i$ coproducts \cite[\S4.5]{benson1998representations}.

\subsection{Monoidal extension and cubical chains}

Let us consider the cellular chains on the interval $\gchains(\gcube)$ as a counital coalgebra in the usual way:
\begin{align*}
\copr[01] &= [0] \ot [01] + [01] \ot [1], &
\copr[0]  &= [0] \ot [0], &
\copr[1]  &= [1] \ot [1] \\
\aug[01] &= 0, &
\aug[0]  &= 1, &
\aug [1]  &= 1.
\end{align*}
This structure can be extended to the chains of cubical sets using the isomorphism
\[
\chains(\cube^n) \cong \gchains(\gcube)^{\ot n}
\]
and the fact that the tensor product of counital coalgebras receives this structure canonically.
Explicitly, for $i \in \{1,2\}$ let $C_i$ be a counital coalgebra, the tensor product $C_1 \ot C_2$ is a counital coalgebra with
\begin{align} \label{e:extension of coproduct}
\copr(c_1 \ot c_2) &= (23) \big( \copr(c_1) \ot \copr(c_2) \, \big), \\
\label{e:extension of augmentation}
\aug(c_1 \ot c_2) &= \aug(c_1) \aug(c_2),
\end{align}
where the symmetric group $\Sym_4$ acts by permuting the tensor factors of $C_1 \ot C_1 \ot C_2 \ot C_2$.

For any cubical set $Y$ the induced structure on its chains agrees with that considered by Serre in \cite{serre1951homologie}, and we refer to it as the \textit{Serre coalgebra} of $Y$.

We can define an $E_\infty$-coalgebra structure extending the Serre coalgebra by describing an extension to all $\chains(\cube^n)$ of the map $\ast \colon \gchains(\gcube)^{\ot 2} \to \gchains(\gcube)$ defined to be non-zero only for
\[
\ast([0] \ot [1]) = [01], \qquad
\ast([1] \ot [0]) = - [01].
\]

For $i \in \{1,2\}$ let $A_i$ be a chain complex equipped with a degree $1$ map $\pr \colon A_i^{\ot 2} \to A_i$ and a chain map $\aug \colon A_i \to \k$ such that $\aug \circ \pr = 0$ and $\partial \pr = \aug \ot \, \id - \id \ot \aug$.
The tensor product $A_1 \ot A_2$ has the same structure, explicitly defined by \eqref{e:extension of augmentation} and
\begin{equation} \label{e:extension of product}
\pr \big( (a_1 \ot a_2) \ot (a_1^\prime \ot a_2^\prime)\big) =
(\id \ot \aug \ot \ast \, + \, \pr \ot \aug \ot \, \id)(23)
\big( a_1 \ot a_2 \ot a_1^\prime \ot a_2^\prime \big)
\end{equation}
where the right-hand side can be given more explicitly by
\[
(-1)^{\bars{a_2} \bars{a_1^\prime}} \big( \aug(a_1^\prime) \, a_1 \ot \ast (a_2 \otimes a_2^\prime)  \, + \, \pr(a_1 \otimes a_1^\prime) \ot \aug(a_2) \, a_2^\prime \big).
\]

Together, formulae \eqref{e:extension of coproduct}, \eqref{e:extension of augmentation}, and \eqref{e:extension of product} induce on the tensor product of $\Med$-bialgebras the same structure.
In particular, $\chains(\cube^n) \cong \gchains(\gcube)^{\otimes n}$ is equipped with an $\Med$-bialgebra structure induced from that in the cellular chain on the interval.
From it, a standard categorical construction -- a Kan extension along the Yoneda embedding~-- provides the chains on any cubical set with the structure of an $E_\infty$-coalgebra extending the Serre coalgebra structure or, more specifically, the structure of a $\UMed$-coalgebra.

Using the monoidal structure on $\Med$-bialgebras we also have a natural $\Med$-bialgebra structure on any tensor product $\chains(\simplex^{n_1}) \otimes \dots \otimes \chains(\simplex^{n_k})$ which defines a natural $\UMed$-coalgebra on multisimplicial chains \cite{gugenheim1957supercomplexes, salvatore2020multisimplicial}

\subsection{Chain contractions for $\Med$} \label{ss:homology of M}

We now explicitly describe chain contractions
\[
\begin{tikzcd}
\arrow[loop, distance=2em, in=195, out=165, "\sigma"'] \Med(s,r) \arrow[r, shift left=3pt, "\proj"] &
\Med(s,r-1) \arrow[l, shift left=3pt, "\incl"]
\end{tikzcd}
\]
for every $s \geq 1$ and $r \geq 0$.
In the next section we will use these to define cup-$(r, i)$ coproducts on $\UMed$-coalgebras.

Recall that $\Med$ is the prop generated by $\copr$, $\aug$, and $\pr$ modulo certain relations (\cref{ss:e-infty generalization}).
Let $\eta \in \Med(s,r)$ be a composition of generators and define
\begin{align*}
\incl(\eta) & = (\id \ot \eta) \circ \copr, \\
\proj(\eta) & = (\aug \ot \, \id^{\ot r-1}) \circ \eta, \\
\sigma(\eta) & = (\pr \ot \, \id^{\ot r-1}) \circ (\id \otimes \eta) \circ \copr.
\end{align*}
These define a chain contraction as above, i.e. they satisfy
\[
\proj \circ \incl - \id = 0
\qquad
\incl \circ \proj - \id = \bd \sigma.
\]

Using the relations defining $\Med$ it is not hard to see that $\Med(s,0) \cong \k$, so we have an explicit chain contraction
\[
\begin{tikzcd}
\arrow[loop, distance=2em, in=195, out=165, "h"'] \Med(s,r) \arrow[r, shift left=3pt, "\proj^r"] &
\k \arrow[l, shift left=3pt, "\incl^r"]
\end{tikzcd}
\]
where
\begin{equation} \label{e:big chain contraction}
h = \sigma + \incl \sigma \proj + \dots + \incl^{r} \sigma \proj^{r}.
\end{equation}
In particular, given that by construction $\UMed(r)$ is a free $\Z[\Sym_r]$-module, these chain contractions shows that $\UMed$ is an $E_\infty$-operad.

\subsection{Steenrod cup-$(r,i)$ products} \label{ss:higher cup-i coproducts}

To generalize the notion of cup-$i$ coproduct structure, consider the cyclic group of order $r$ and the minimal free resolution of $\Z$ as a $\Z[\Cyc_r]$-module
\begin{equation} \label{e:minimal resolution r>2}
\begin{tikzcd} [column sep = .5cm]
\mathcal W(r) = \Z[\Cyc_r]\{e_0\} & \arrow[l, "\,T"'] \Z[\Cyc_r]\{e_1\} & \arrow[l, "\,N"'] \Z[\Cyc_r]\{e_2\} & \arrow[l, "\,T"'] \cdots
\end{tikzcd}
\end{equation}
where
\begin{equation} \label{e:transposition and norm maps}
N = 1 + \rho + \dots + \rho^{r-1}, \qquad
T = \rho - 1,
\end{equation}
and $\rho$ is a generator of $\Cyc_r$.

Let $C$ be a $\UMed$-coalgebra, for example the chains on a simplicial or cubical set.
For $r \geq 2$, let $\psi(r) \colon \cW(r) \to \End(C, C^{\ot r})$ be the $\Cyc_r$-equivariant chain map defined recursively by
\begin{equation}\label{e:definition of psi}
\begin{split}
\psi(r)(e_0) & = \Delta^{r-1}, \\
\psi(r)(e_{2m+1}) & = h\,T\,\psi(r)(e_{2m}), \\
\psi(r)(e_{2m}) & = h\,N\,\psi(r)(e_{2m-1}),
\end{split}
\end{equation}
where $T$ and $N$ are explicitly defined in \eqref{e:transposition and norm maps} and $h$ in \eqref{e:big chain contraction}.
The \textit{Steenrod cup-}$(r, i)$ \textit{product} of $C$ is defined for every $r, i \geq 0$ as the image in $\End(C, C^{\ot r})$ of $\psi(e_i)$.
We remark that for $r = 2$ the resulting cup-$i$ coproduct structure on simplicial chains recovers Steenrod's original construction \eqref{e:cup-i coproducts}.

\subsection{Steenrod operations} \label{ss:steenrod operations}

Let $p$ be an odd prime.
We now review a construction of Steenrod operations in mod $p$ cohomology analogue to the one given in \cref{ss:steenrod squares} for Steenrod squares.
We remark that Steenrod square operations are parameterized by the mod $2$ homology of $\Cyc_2 = \Sym_2$.
As explained for example in \cite[Corollary~VI.1.4]{adem2004milgram}, an inclusion $\Cyc_p \to \Sym_p$ induces a surjection in mod $p$ homology.
Using the homological degrees where this surjection is non-zero we have the following construction explained in more detain in \cite{steenrod1952reduced, steenrod1953cyclic, steenrod1962cohomology, may1970general}.

Let $C$ be a $\UMed$-coalgebra and $C^\vee = \Hom(C, \Fp)$.
We simplify notation and denote the explicit linear map $\psi(p)(e_i) \colon C \to C^{\ot p}$ defined in \eqref{e:definition of psi} simply by $\psi_i$.
For any integer $s$, the \textit{Steenrod operations}
\begin{equation*}
P_s \colon \rH(C^\vee; \Fp) \to \rH(C^\vee; \Fp)
\end{equation*}
and
\begin{equation*}
\beta P_s \colon \rH(C^\vee; \Fp) \to \rH(C^\vee; \Fp)
\end{equation*}
are defined by sending the class represented by $\alpha \in C^\vee$ of degree $q$ to the classes represented respectively for $\varepsilon \in \{0,1\}$ by
\begin{equation*}
\pm \, \alpha^{\ot p} \circ \psi_{(2s-q)(p-1)-\varepsilon}(-)
\end{equation*}
where a possible sign convention is introduced and motivated in \cite[(6.1)]{steenrod1953cyclic}.

\subsection{Cartan and Adem relations}

There is a conceptual pattern producing additional homological information from relations bounding primary structure.
As an example, we have seen that Steenrod operations can be understood as structure on cohomology deduced from lifting to the cochain level the commutativity relation satisfied by the cup product.
Steenrod operations, which we now regard as primary cohomological structure, satisfy relations that lead to finer structure.
The first of these is the \textit{Cartan relation} that establishes a connection between Steenrod operations and the cup product, it is given by:
\begin{equation*}
P_s\big(\alpha \beta \big) = \sum_{i+j=s} P_i\big( \alpha \big) P_j\big( \beta \big).
\end{equation*}
The second is the \textit{Adem relation}, controlling the iteration of Steenrod operations, and given by:
\begin{enumerate}
	\item If $p = 2$ and $a > 2b$, then
	\begin{equation*}
	P_{a} P_{b} = \sum_i \binom{2i-a}{a-b-i-1} P_{a+b-i}P_i,
	\end{equation*}
	\item If $p > 2$ and $a > pb$, then
	\begin{equation*}
	P_{a} P_{b} = \sum_i (-1)^{a+i} \binom{pi-a}{a-(p-1)b-i-1} P_{a+b-i}P_i.
	\end{equation*}
\end{enumerate}
There are version of these using the $\beta P_s$ operations but we do not write them here, see for example \cite{steenrod1962cohomology, may1970general}.

Steenrod operations and the above relations are a homological consequence of an $E_\infty$-structure.
For cellular chains we have seen explicit cochain level constructions, the cup-$(p,i)$ products, inducing the Steenrod operations, and it is desirable to produce cochains enforcing these relations.
For the even prime case, Cartan and Adem coboundaries have been constructed effectively in \cite{medina2020cartan} and \cite{medina2021adem} respectively.
Cartan coboundaries for odd primes can be constructed with the tools already described, but the Adem relation requires additional techniques not yet available.

A source of motivation for these cochain level constructions comes from their use in the study of topological phases as we overview next.

\subsection{Symmetry protected topological phases and cochain constructions}

A central problem in physics is to define and understand the moduli ``space'' of quantum systems with a fixed set of invariants, for example their dimension and symmetry type.
In condensed matter physics, quantum systems are
presented using \textit{lattice models} which, intuitively, are given by a Hamiltonian presented as a sum of local terms on a Hilbert space associated to a lattice in $\R^n$.
We think of these as defined on flat space.
One such system is said to be \textit{gapped} if the spectrum of the Hamiltonian is bounded away from $0$, and two Hamiltonians represent the same \textit{phase} if there exists a deformation between them consisting only of systems that remain bounded from below.

Given a lattice model, using cellular decompositions and state sum type constructions one can often compute the corresponding \textit{partition functions} on spacetime manifolds from actions expressed in terms of gauge fields represented by cochains and cochain level structures: Stiefel--Whitney cochains, cup-$i$ products and Cartan/Adem coboundaries, for example.
Subdivision invariance gives rise to a topological quantum field theory, which in the \textit{invertible} case is expected to be controlled by a generalized cohomology theory \cite{xiong2018minimalist, gaiotto2019cohomology, freed2021reflection}.
The cochain level structure used in the definition of the cellular gauge theory is interpreted from this point of view as describing a cochain model of the Postnikov tower of the relevant spectrum.
For example, fermionic phases protected by a $G$-symmetry are believed to be classified by applying to $BG$ the Pontryagin dual of spin bordism \cite{kapustin2015cobordism, kapustin2017fermionic}.
Building on these insights, A. Kapustin proposed a structural ansatz in low dimensions that Greg Brumfiel and John Morgan verified by constructing cochain models of certain connective covers of said spectrum \cite{brumfiel2016pontrjagin, brumfiel2018pontrjagin}.
The resulting models use cup-$i$ products and Adem coboundaries to represent $k$-invariants and the additive structure of the spectrum.

	\sloppy
	\printbibliography
\end{document}